\documentclass[12pt]{article}
\usepackage{amssymb,latexsym,amsmath,mathrsfs,graphicx}
\begin{document}
\title{Concerning life annuities\footnote{Originally published as
\emph{Sur les rentes viageres}, Memoires de l'academie des sciences de
Berlin \textbf{16} (1767), 165-175. A copy of the original
text is available electronically at the
Euler Archive, at www.eulerarchive.org. This paper is E335 in the Enestr\"om 
index.}}
\author{Leonhard Euler\footnote{Date of translation: February 17, 2005.
Translated from the French by
Christian L\'eger, 3rd year undergraduate in Honours Mathematics,
and Jordan Bell, 2nd year undergraduate in Honours Mathematics,
School of Mathematics and Statistics, Carleton University, Ottawa, Ontario,
Canada. Email: cleger@connect.carleton.ca and jbell3@connect.carleton.ca }}
\date{}
\maketitle

1. Having established the true principle for which one must
base the calculation for life annuities, I believe that the
development of this calculation will not fail to be very
interesting, indeed as much for those who wish to run such an establishment
as for those who will wish to profit from it. I worked on this matter
in my ``Recherches g\'en\'erales sur la mortalit\'e et la multiplication
du genre humain'', where I showed the correct
method of determining
by calculation how much a man of a given age must pay in order to enjoy
for his whole life a specific annual income. However, since this calculation
seemed to me at the time very cumbersome, I could not convince myself
to carry it out. Then, on a given occassion I was forced to undertake this
calculation, for which by means of certain tricks to abbreviate the calculation,
I was fortunately able to complete it. 

2. There are two things on which the calculation of the income from
life annuities must be based: one of them is good data on mortality,
which tells us for each age how many will probably die during the course of
one or many years; the other is the way in which the manager can make
the money grow in value that he will have received from the annuitants,
or at which interest rates he is able to place it. Together, these two
things essentially determine what returns the manager will be able to commit himself
to, as much with respect to the amount initially deposited as with respect to the
age of the annuitant, because it is evident that the more the profit
the manager can obtain from the capital in his hands, the more they will
also be able to provide the annuitant strong returns.

3. For the mortality list, the manager would no doubt risk much
if they based it on the mortality of men in general
which we compile from observations made in a big city or in a whole country,
where we keep track equally of all men, vigorous and infirm. Now, when
the business is to procure oneself life annuities, it is very natural
that necessarily exluded are all those whose constitution does not seem to
promise long life, thus there is reason to regard annuitants as a more
robust species. It is also with this consideration that I chose
in my ``M\'emoire allegu\'e la liste de M. Kerseboom'', which he
generated from observations made only on persons who enjoy
life annuities: and using also this same list will assist me in
the elaboration of the following calculations.

4. If the manager is not well able to place the capital that is given
him by the annuitants, he would be able to afford naught very mediocre
returns such that no one would want them. Once the city of Amsterdam
paid ten percent annuities to all persons under the age of twenty, that
is
for 1000 florins they are paid 100 per year, which is an annuity
so rich that the city would have suffered considerable loss if it had
not generated almost 10 percent per year of the fund which this enterprise
had provided it. Thus, if we can only count on 5 percent interest, the
annuities must become far less considerable; however, it is on this
which it seems annuities must be fixed, provided that those
who will have occassion to obtain greater profit from it will hardly
participate in such an enterprise, which could only be worthwhile
after a large number of years.

5. For determining the cost of these annuities,
we fix for each age an average lifespan, which one is equally
likely to survive as to die before having reached it; that is, this term is chosen
such that as many men of the same age die before this term as die after.
Thus we presume that all men of this age reach exactly this term, and that
they then die; on this we believe we can fix with certainty the price
of the annuities, since the value of the annuities must be payable during
a given number of consecutive years: and we estimate that the profit
that the manager obtains on the side of those who die before their
projected term is exactly compensated by the loss caused by those
annuitants who survive this term. However, we will easily understand that
this reasoning is flawed, since it does not take into account the
reduction in current price of an annuity which will not be paid out before
many years. Given this circumstance, it will be necessary to base
the calculation on true principles, as I explained in my
aforementioned ``M\'emoire'', without making use of any reasoning
which might seem suspect.

6. To achieve this, we consider a number of 1000 children born at the same
time,  and that the symbols (1), (2), (3), (4), etc. indicate the number
of those who still live at the end of 1, 2, 3, 4, etc. years, 
such that in general $(m)$ represents the number of those who will
obtain the age of $m$ years. Now let $r$ be the annual annuity that a man
of $m$ years wishes to receive, and let $x$ be the price which he must then
pay at the present to the manager, which must be a just equivalent
of the spending to which the manager engages himself by this agreement.
To determine this price $x$, many men of the same age $m$ must be considered,
and those who reach this age. Let $(m)$ be this number of men, and the sum
they will presently pay to the manager will be equal to $(m)x$, which 
must be sufficient to provide for all the annuities which he will have to pay
afterwards.

7. Of these $(m)$ men, there will remain alive after one year $(m+1)$, after
two years $(m+2$), after three years $(m+3)$, and so on. Thus the manager
will have to pay after one year $(m+1)r$, after two years ($m+2)r$, after three
years $(m+3)r$, etc., until all of these annuitants will be extinguished.
We thus only have to reduce all of the payments at the present time by the amount
of 5 percent, and make the sum equal to $(m)x$ to determine the value of $x$.
Now to make the calculation more general, instead of $\frac{105}{100}$ or
$\frac{21}{20}$, let us write the letter $\lambda$, and the sum of all
the annuities which the manager must pay successively will now be:
\[ \frac{(m+1)r}{\lambda}+\frac{(m+2)r}{\lambda^2}+\frac{(m+3)r}{\lambda^3}+
\frac{(m+4)r}{\lambda^4}+\textrm{ etc.} \]
which being equal to $(m)x$, will give:
\[ x=\frac{r}{(m)} \Bigg ( \frac{(m+1)}{\lambda}+\frac{(m+2)}{\lambda^2}+
\frac{(m+3)}{\lambda^3}+\frac{(m+4)}{\lambda^4}+\textrm{ etc.} \Bigg )\]

8. Then the exact price is found which a man of $m$ years must pay in order
to enjoy an annual annuity $r$ during his whole life, this one having been
initially placed at 5 percent, puts the manager precisely within the means
of paying from that point the annuities as long of the number of annuitants
is sufficiently large. We understand well, having thus placed initially
all the capital which the manager will have received, the following year
the interest will not be sufficient to pay the annuities but that it will
be needed to employ part of the capital, hence the capital will suffer every
year a diminuition: however, it will only be entirely extinguished when the
annuitants are dead. For this reason the manager will be well obliged to
raise the price of the annuities that I have just found, according to
the particular circumstances and expenses which such an 
establishment requires. 

9. We clearly see that the determination of this price called $x$ requires
a calculation as tedious as it is unpleasant, especially for low ages, where
the number of terms to be added together is very considerable. But it is
not hard to notice, that having already done a calculation for a certain
age, we will from it be easily able to extract the one which corresponds
to a later or earlier year. To better explain this artifice, I will employ this
character $\overline{m}r$ to indicate the price which a man of age $m$
must pay for the annuity $r$: in order that
\[ \overline{m}= \frac{1}{(m)} \Bigg ( \frac{(m+1)}{\lambda} +
\frac{(m+2)}{\lambda^2}+\frac{(m+3)}{\lambda^3}+\frac{(m+4)}{\lambda^4} +
\textrm{etc.} \Bigg ), \]
from there, for men aged $m+1$ years we will have,
\[ \overline{m+1}= \frac{1}{(m+1)} \Bigg ( \frac{(m+2)}{\lambda}+ \frac{(m+3)}{\lambda^2} + \frac{(m+4)}{\lambda^3} + \frac{(m+5)}{\lambda^4}+ \textrm{etc.} \Bigg ), \]
from which we conclude: 
\[ \lambda(m)\overline{m} = (m+1)+(m+1)\overline{m+1}, \]
and starting,
\[ \overline{m}=\frac{1}{\lambda} \cdot \frac{(m+1)}{(m)}(1+\overline{m+1}), \]
such that having found the value $\overline{m+1}$, we will from it calculate
easily enough the value of $\overline{m}$.

10. With the aid of this artifice, after having started with the age of 90
years, I calculated the price of an annuity $r$ successively for all inferior ages,
down to those children newly born; where I obtained the following
table, by fixing the annuity $r$ at 100 crowns, and interest at 5 percent.

\begin{center}
Table

{\em Table which indicates the price of a life annuity of 100 crowns
for all ages}
\begin{tabular}{p{1cm}|p{2cm}|p{2cm}||p{1cm}|p{2cm}|p{2cm}}
age in years&number of survivors&price of annuity&age in years
&number of survivors&price of annuity\\
0&1000&1155.50&25&552&1403.60\\
\hline
1&804&1409.04&26&544&1395.45\\
2&768&1448.84&27&535&1389.87\\
3&736&1487.43&28&525&1387.16\\
4&709&1521.27&29&516&1382.54\\
5&690&1541.32&30&507&1376.82\\
\hline
6&676&1551.90&31&499&1368.84\\
7&664&1558.94&32&490&1363.68\\
8&654&1561.92&33&482&1355.63\\
9&646&1560.33&34&475&1344.38\\
10&639&1556.29&35&468&1332.71\\
\hline
11&633&1549.59&36&461&1320.60\\
12&627&1542.64&37&454&1308.01\\
13&621&1535.42&38&446&1298.04\\
14&616&1525.28&39&439&1284.67\\
15&611&1514.65&40&432&1270.76\\
\hline
16&606&1503.50&41&426&1253.09\\
17&601&1491.81&42&420&1234.54\\
18&596&1479.54&43&413&1218.24\\
19&590&1469.31&44&406&1201.21\\
20&584&1458.63&45&400&1180.19\\
\hline
21&577&1450.18&46&393&1161.27\\
22&571&1438.68&47&386&1141.44\\
23&565&1426.66&48&378&1123.88\\
24&559&1414.07&49&370&1105.59\\
25&552&1403.60&50&362&1086.52\\
\hline
$m$&$(m)$&$\overline{m}$&$m$&$(m)$&$\overline{m}$
\end{tabular}
\end{center}

\begin{tabular}{p{1cm}|p{2cm}|p{2cm}||p{1cm}|p{2cm}|p{2cm}}
age in years&number of survivors&price of annuity&age in years
&number of survivors&price of annuity\\
50&362&1086.52&70&175&638.30\\
\hline
51&354&1066.62&71&165&610.83\\
52&345&1049.17&72&155&582.75\\
53&336&1031.14&73&145&554.09\\
54&327&1012.49&74&135&524.89\\
55&319&989.78&75&125&495.22\\
\hline
56&310&969.44&76&114&470.16\\
57&301&948.35&77&104&441.13\\
58&291&929.98&78&93&417.98\\
59&282&907.64&79&82&397.75\\
60&273&884.44&80&72&375.64\\
\hline
61&264&860.32&81&63&350.77\\
62&254&838.90&82&54&329.69\\
63&245&813.21&83&46&309.38\\
64&235&790.20&84&39&279.44\\
65&225&766.59&85&32&257.60\\
\hline
66&215&742.30&86&26&232.90\\
67&205&717.43&87&20&217.91\\
68&195&691.93&88&15&205.07\\
69&185&665.14&89&11&193.62\\
70&175&638.30&90&8&179.54\\
\hline
$m$&$(m)$&$\overline{m}$&$m$&$(m)$&$\overline{m}$
\end{tabular}

11. M. Kerseboom only continued the table on mortality up to 95 years,
and for this reason I did not judge it convenient to continue this one
beyond 90 years, since probably at this age nobody will desire
life annuities. At the least, in almost all plans, such ancients find themselves
filed into the same class as those of 60 or 70 years, notwithstanding
that it would be very unjust if we wanted to demand of a
nonagenarian more than a third of the price which a septuagenarian must
pay or a quarter of what a sexagenarian must pay. However, if we are curious
to see the continuation of my table, here it is:

\begin{tabular}{l|l|l|l|l|l}
$m$&90&91&92&93&94\\
$(m)$&8&6&4&3&2\\
$\overline{m}$&179.54&151.35&138.38&93.73&47.62
\end{tabular}

But I would not advise a manager to get involved with such 
ancients unless their number be sufficiently considerable; which is a general
rule for all establishments founded on probabilities.

12. From there we will conclude easily how much the manager should pay
in interest for each age for a given sum which we will have initially
given. It is not necessary to enter here in as much detail, and it suffices
to mark for every five years which the annuitants might expect.

\begin{tabular}{r|r||r|r||r|r}
age&percent&age&percent&age&percent\\
0&$8 \frac{2}{3}$&30&$7 \frac{1}{4}$&60&$11 \frac{1}{3}$\\
5&$6 \frac{1}{2}$&35&$7 \frac{1}{2}$&65&13\\
10&$6 \frac{1}{3}$&40&8&70&$15 \frac{2}{3}$\\
15&$6 \frac{1}{2}$&45&$8 \frac{1}{2}$&75&20\\
20&$6 \frac{3}{4}$&50&9&80&$25 \frac{2}{3}$\\
25&$7$&55&10&85&$38 \frac{1}{2}$\\
30&$7 \frac{1}{4}$&60&$11 \frac{1}{3}$&90&$55 \frac{1}{2}$
\end{tabular}

On this basis, the manager will obtain no profit unless he is able
to appreciate his money at more than 5 percent.

13. Thus, if a state has need of money and it can find at better than
5 percent interest as much as it needs, then it would be assuredly
very bad if it wished to establish such life annuities that I have
just determined on this basis of 5 percent, since, with regard to
the burden which such an establishment necessarily imposes, it would
always be better to borrow the sum which it needs at 5 percent which it
would then be able to satisfy  according to circumstance, in the place
of life annuities which would remain in its charge during a very long
period. Or else, the price on the annuities would have to be raised
beyond that which I have just described to bring it any benefit, but then
it might be difficult to find any more annuitants unless it be 
ancients over 60 years old, who could be amazed by
 interests of 10 percent or more.

14. But wanting to establish more advantageous life annuities for 
annuitants, it would be a project hardly proper for satisfying a state,
since this would return to the same, then if one wanted to burden
one's self with debts at 6 or more percent: while we could borrow
at 5 percent without being subjected to the trouble which life
annuities would incur. Indeed, if a state wants to establish
the life annuities described here and calculuated on the basis of 5
percent, this charge need only be regarded as a loan taken at 6 percent,
because of the number of arrangements which would be required for it.
Thus I hardly see a case where the establishment of life annuities
can be advantageous to a state , as long as we can borrow money at 
5 percent and perhaps less. But we can imagine a different
type annuity which may be better tailored, although still based
on 5 percent. I want to speak of annuities which must not
begin to run before 10 or 20 years; and we easily understand that the
price of such annuities will be quite mediocre, capable of
attracting the public only poorly.

15. Consider this question also in general, and ask how much a man of age
$m$ years must pay presently to obtain an annual annuity $r$ which will
not begin to be paid to him before $n$ years, such that after this time
he is able to enjoy it regularly until his death. Let $x$ be the current
price of this annuity, and we will find as here below:
\[
x=\frac{r}{(m)} \Big ( \frac{(m+n)}{\lambda^n}+\frac{(m+n+1)}{\lambda^{n+1}}+
\frac{(m+n+2)}{\lambda^{n+2}}+\textrm{ etc.} \Big ).
\]
Thus to calculate the previous explained ordinary annuities, we will
have:
\[
\overline{m+n-1}=\frac{1}{(m+n+1)} \Big ( \frac{(m+n)}{\lambda}+
\frac{(m+n+1)}{\lambda^2}+\frac{(m+n+2)}{\lambda^3}+\textrm{ etc.} \Big ),
\]
from which we conclude,
\[
x=\frac{r}{(m)}\cdot \frac{(m+n-1)}{\lambda^{n-1}}\cdot \overline{m+n-1}
=
\frac{r}{\lambda^{n-1}}\cdot \frac{(m+n-1)}{(m)}\cdot \overline{m+n-1},
\]
where $\overline{m+n-1}r$ expresses the current price of the ordinary
annuity for a man of age $m+n-1$ years.

16. Thus if we ask the current price of an annual annuity of 100 crowns
which will not begin to be paid for 10 years, for a man of age $m$ years
we will take from the table developed in  \S 10 the price of the ordinary
annuity which is suited to the age $m+9$ years, and 
we will multiply by $\Big ( \frac{20}{21} \Big )^9 \frac{(m+9)}{(m)}$
to obtain the value sought for $x$. From there I calculated the following
tables for every 5 years:
\newpage

\begin{center}
Table

{\em Prices of life annuities of 100 crowns which will not start
to run for 10 years}
\begin{tabular}{p{1cm}|p{2.5cm}||p{1cm}|p{2.5cm}||p{1cm}|p{2.5cm}}
age years&price of the annnuity&age year&price of the annuity&age year&price of the annuity\\
0&649.75&30&717.05&60&290.55\\
5&877.77&35&671.73&65&203.11\\
10&874.50&40&610.40&70&120.14\\
15&833.95&45&533.55&75&56.20\\
20&787.43&50&455.78&80&19.07\\
25&745.72&55&375.25&&\\
30&717.05&60&290.55&&
\end{tabular}
\end{center}

\begin{center}
Table

{\em Price of a life annuity of 100 crowns will not begin
to run for 20 years}
\begin{tabular}{p{1cm}|p{2.5cm}||p{1cm}|p{2.5cm}||p{1cm}|p{2.5cm}}
age year&price of the annuity&age year&price of the annuity&age year&price of the annuity\\
0&343.06&30&319.30&60&47.28\\
5&453.36&35&272.96&65&19.17\\
10&441.81&40&234.47&70&4.82\\
15&413.60&45&183.72&&\\
20&382.17&50&134.52&&\\
25&349.63&55&87.91&&\\
30&319.30&60&47.28&&
\end{tabular}
\end{center}

17. Maybe such a project of life annuities will better succeed 
in spite of being fixed on a basis of 5 percent. It seems that it will
always be advantageous for a newborn child to ensure it, assuming a price
of 343 or 350 crowns, a fixed annuity of 100 crowns per year, though it
will not start to be paid before the child attains the age of 20:
and if we wish to use the sum of 3500 crowns, it will always be
a nice establishment to be able to enjoy from the age of 20 years a fixed
pension of 1000 crowns. However, it is doubtful that there will be many
parents who will well be willing to make such a sacrifice for the good of their
children. Perhaps there will be more men of 60 years who will not mind
paying first 3000 crowns to be assured a fixed pension of 1000 crowns per year
once they have passed their 70-th year. 

\end{document}